\pdfoutput=1
\documentclass[10pt]{article}
\usepackage{times,amsmath,amssymb,exscale,mathrsfs}
\usepackage{graphicx}
\usepackage{xcolor} 
\usepackage{epsfig}
\newcommand{\C}{\mathbb{C}}

\newcommand{\R}{\mathbb{R}}

\newcommand{\W}{\mathcal{W}}

\title{Vibrational resonance: a study with high-order word-series averaging}
\author{A. Murua\footnote{Konputazio Zientziak eta A.\ A.\  Saila, Informatika
 Fakultatea, UPV/EHU, E--20018 Donostia--San Sebasti\'{a}n,  Spain. Email: Ander.Murua@ehu.es}
\  and J.M. Sanz-Serna\footnote{(Corresponding author) Departamento de  Matem\'aticas, Universidad Carlos III de Madrid, Avenida de la Universidad 30, E-28911Legan\'es (Madrid), Spain.
 Email: jmsanzserna@gmail.com}}
\date{\today}
\begin{document}
\maketitle

\begin{abstract}
We study a model problem describing vibrational resonance  by means of
a high-order averaging technique based on so-called word series. With the technique applied here, the tasks of constructing the averaged system and  the associated change of variables  are divided into two parts. It is first necessary to build recursively a set of so-called {\em word basis} functions and, after that, all the required manipulations involve only scalar coefficients that are computed by means of simple recursions. As distinct from the situation with other approaches, with word-series, high-order averaged systems may be derived without having to compute the associated change of variables. In the system considered here, the construction of high-order averaged systems  makes it possible to obtain very  precise approximations to the true dynamics.
\end{abstract}

\medskip

\noindent{\bf Keywords and sentences:} High-order averaging, vibrational resonance, formal series
\medskip


\section{Introduction}
We study a model problem describing vibrational resonance \cite{vr} by means of
the high-order averaging technique introduced in the series of articles \cite{part1}, \cite{part2}, \cite{orlando}, \cite{part3} (see \cite{kurusch} for a summary).

In  devices with vibrational resonance \cite{vr}, \cite{baltanas}, \cite{raj} the response of a system driven by a low-frequency forcing may be enhanced by the presence of high-frequency vibrations of suitable amplitude. Such devices feature  in several current applications, including energy harvesting \cite{eh}, \cite{eh2}, where the aim is to exploit the energy contained in the vibrations by converting it into electrical energy. The corresponding systems of differential equations  are {\em highly oscillatory} in the sense that their solutions include periods much shorter than the time interval of interest. The numerical integration of such systems may be a very expensive task, as, typically, numerical integrators have to operate with step sizes significantly shorter than the shortest period present in the solution and the simulation has to be carried out over many periods. The literature contains many suggestions of nonstandard integrators to be applied in highly oscillatory situations, including the heterogeneous multiscale method \cite{etsai}, \cite{ima}, methods with multiple time steps \cite{sinum}, multirevolution techniques \cite{sam}, etc. However the use of such ad hoc integrators is prone to unexpected instabilities and inaccuracies \cite{ii}.

Averaging techniques are very useful to deal with  oscillatory problems; they may be applied in combination with numerical integration, to \lq soften\rq\ the difficulties of the system to be simulated numerically, or in a purely analytical way. There are many techniques to carry out the computations required by the method of averaging (see e.g.\ \cite{averbook}); here we use the {\em word series} \cite{china} approach introduced in \cite{part1}, \cite{part2}, \cite{orlando}, \cite{part3}. With the technique applied here, the tasks of constructing the averaged system and  the associated change of variables  are divided into two parts. It is first necessary to build recursively a set of so-called {\em word basis} functions and, after that, all the required manipulations involve only scalar coefficients that are computed by means of simple recursions. As distinct from the situation with other approaches, with word-series high-order averaged systems may be derived without having to compute the associated change of variables.
The use of word series techniques is not confined to averaging; word series may be applied to compute normal forms of discrete or continuous dynamical systems \cite{words}, \cite{juanluis}, \cite{kurusch}, find formal invariants of differential systems \cite{juanluis},  analyze numerical integrators \cite{words}, \cite{alfonso}, etc.

Section 2 provides a brief summary of the averaging technique we employ. It should be emphasized that while, for the sake of simplicity, this summary is restricted to the case of {\em periodic} forcing, the technique may be applied to {\em quasiperiodic} problems without essential changes \cite{part2}, \cite{kurusch}. Section 3 studies a model problem taken from the original reference on vibrational resonance by Landa and McClintock \cite{vr}. We show how averaging provides insight into the mechanism causing vibrational resonance; this mechanism is related to the creation of an effective potential, akin e.g. to that responsible for stabilizing Kapitza's inverted pendulum \cite{ima}.
Furthermore, the construction of high-order averaged systems  makes it possible to get very  precise approximations to the true dynamics of the system under investigation.

The example here is based on the overdamped oscillator; the technique may be equally applied to underdamped models \cite{eh}, \cite{eh2}. Work to extend our study to realistic devices that exhibit vibrational resonance is under way.

\section{Averaging with word series}
After Fourier expansion, we may assume that  the  periodic problem  to be studied is of the form:
\begin{equation}\label{eq:edo}
\frac{d}{dt} x = \sum_{k=-\infty}^\infty \exp(ik\omega t)f_k(x), \qquad x(t_0) = x_0\in\R^D.
\end{equation}
In the averaging technique used here,
the elements (indices) \( k = 0,\pm 1,\pm 2,\dots\) in \eqref{eq:edo} are seen as the {\em letters} of an alphabet and each (possibly empty) string
of letters \( k_1k_2\dots k_n\), \( n = 0, 1,\dots\), is called a {\em word}. The symbol \(\W\) is used to denote the set of all words. With each word \(w\in\W\) we associate a {\em word basis function} \(f_w\). For the empty word, \(f_\emptyset(x)\) is the identity map \( x\mapsto x\) in \(\R^D\). For nonempty words the \(f_w(x)\) are constructed recursively from the \( f_k(x)\) that feature in \eqref{eq:edo}. The recipe is
\[
f_{k_1\dots k_n}(x) = f^\prime_{k_2\dots k_n}(x)f_{k_1}(x),
\]
where \(f^\prime_{k_2\dots k_n}(x)\) is the Jacobian matrix of \(f_{k_2\dots k_n}(x)\). We  define \( \C^\W\) as the set (vector space) of all mappings \( \delta : \W\rightarrow \C\); if \(\delta\in\C^\W\) and \(w\in\W\), \(\delta_w\) represents the complex number that \(\delta\) associates with \(w\). To each \(\delta\in\C^\W\) there corresponds
a {\em word series}
\[
W_\delta(x) = \sum_{w\in\W} \delta_wf_w(x).
\]

With this terminology, the main result from \cite{part2} is as follows. Consider the {\em averaged} problem
\begin{equation}\label{eq:edoaver}
\frac{d}{dt} X = W_{\bar \beta(t_0)}(X),\qquad X(t_0) = x_0
\end{equation}
and the time-dependent {\em change of variables}
\begin{equation}\label{eq:change}
 x = W_{\kappa(t\omega,t_0)}(X),
\end{equation}
where the coefficients \(\bar\beta_w(t_0)\), \(\kappa_w(t\omega,t_0)\), \(w\in\W\), necessary to write the word series may be computed explicitly (see below). Then the solution \(x(t)\) of \eqref{eq:edo} may be represented as
\begin{equation}\label{eq:main}
x(t) = W_{\kappa(t\omega,t_0)}(X(t)),
\end{equation}
where \( X(t)\) is the solution of \eqref{eq:edoaver}. The change of variables \eqref{eq:change} is \( 2\pi/\omega\) periodic in \(t\) and, in addition, at the stroboscopic times \( t_\ell = t_0+\ell(2\pi/\omega)\),
\(x\) and \( X \) coincide. In general the series in the right hand-sides of \eqref{eq:edoaver} and \eqref{eq:change} do not converge and have to be truncated as demonstrated in the next section. With the help of such truncations it is possible to approximate \( x(t)\) accurately without solving the oscillatory problem \eqref{eq:edo}.

The coefficients \(\bar\beta_w(t_0)\) are computed recursively by means of the formulas
\begin{eqnarray*}
\bar \beta_{k}(t_0)&=& 0, \\
\bar \beta_{ 0}(t_0) &=& 1, \\
\bar \beta_{ 0^{r+1}}(t_0) &=& 0, \\
\bar \beta_{{0}^r k}(t_0)
&=&
\displaystyle\frac{i}{k \omega}( \bar
\beta_{ { 0}^{r-1}k}(t_0) -  \bar \beta_{{ 0}^{r}}(t_0)e^{ik \omega t_0}),\\
\bar \beta_{k \ell_1\cdots \ell_s}(t_0)
&=&
\displaystyle\frac{i}{k \omega}(e^{i k \omega t_0}
\bar \beta_{\ell_1\cdots \ell_s}(t_0)- \bar \beta_{(k+\ell_1) \ell_2\cdots \ell_s}(t_0)),\\
\bar \beta_{{ 0}^r k \ell_1\cdots \ell_s}(t_0)
&=&
 \displaystyle\frac{i}{k \omega}(
\bar \beta_{{ 0}^{r-1} k \ell_1\cdots \ell_s}(t_0) -  \bar \beta_{{ 0}^{r} (k+\ell_1) \ell_2\cdots \ell_s}(t_0)).\\
\end{eqnarray*}
Here the integer \(k\) is \(\neq 0\) and \(0^r\), \(r>0\), denotes the word consisting of \( r\) zeros. Note that, from these formulas, \( \bar\beta_w\) is of size \(\mathcal{O}(1/\omega^{n-1})\) for  \( n\) letter words \(w\), \( n >0\). The, very similar formulas for the coefficients \(\kappa(t\omega,t_0)\)
may be seen in \cite{kurusch}.

\section{Application to a vibrational resonance  problem}
The following  one-dimensional overdamped double-well oscillator
\begin{equation}\label{eq:vr1}
\frac{d z}{dt} = z -z^3 +A \cos \nu t + B \omega \cos \omega t, \qquad z(0) = z_0,
\end{equation}
has been considered in \cite{vr} as a simple model to demonstrate vibrational resonance. Here \( 0 <\nu\ll \omega\), the parameter \( A \) represents the amplitude of the applied driving force and \( B\) measures the size of the fast background vibration. We note that in \cite{vr} the vibrational term is written in the alternative format \( C \cos \omega t\) rather than as
\(B \omega \cos \omega t\). Since no hypotheses are made in what follows as to the size of  \(B\), both formats are equivalent; when \(C\) and \(\omega\) are given specific numerical values the problem  considered in \cite{vr} may be cast in the form \eqref{eq:vr1} by defining \( B = C/\omega\).
 However writing the amplitude of the vibration as
\(B \omega\) is more meaningful than writing it as \(C\) because if one sees \(\omega\) as a parameter and keeps \( C\) constant then it is clear that for \(\omega\) sufficiently large the vibrational term \( C\cos \omega t\) cannot be expected to exert any significant influence as it converges weakly to zero.

We begin by removing from \eqref{eq:vr1} the \(\mathcal{O}(\omega)\)  term by means of the preliminary change of variables
\(
z = x + B \sin \omega t
\)
(\(z\) and \(y\) coincide at the stroboscopic times \(\ell(2\pi/\omega)\)), which leads to
\begin{eqnarray}\label{eq:vr2}
\frac{d y}{dt} &= &y -\frac{3}{2} B^2 y  -y^3+A \cos \nu t
\\
\nonumber
&&+ B \left( 1-\frac{3}{4} B^2-3y^2\right)\sin \omega t+
\frac{3}{2}  B^2 y\cos 2\omega t +\frac{1}{4} B^3 \sin 3\omega t.
\end{eqnarray}
Note that the conservative force \( z-z^3\) in \eqref{eq:vr1} has been \lq softened\rq\ to \((1-3B^2/2)y -y^3\) in \eqref{eq:vr2}. In terms of the corresponding potentials, we have moved from \( -z^2/2+ z^4/4\), with wells of depth
\(1/4\) at \( z=\pm 1\) to the new symmetric {\em effective} potential
\[
\left( -\frac{1}{2} +\frac{3}{4}B^2\right) y^2 +\frac{1}{4} y^4.
\]
As \(|B|\) increases from \(B = 0\), the depth of the wells decreases and for \(|B|\geq \sqrt{2/3}\) the two potential minima are merged into  a single  minimum at \( y = 0\).

By considering an auxiliary variable \(\phi\) with \((d/dt)\phi= \nu\) and the vector \( x = (\phi,y)^T\in\R^2\), \eqref{eq:vr2} takes the form \eqref{eq:edo} studied in the preceding section with
\begin{eqnarray*}
f_0(x) & =& (\nu, y  -\frac{3}{2} B^2 y-y^3 +A \cos \phi)^T,\\
f_1(x) & = & (0, -\frac{i}{2}B +\frac{3i}{8}B^3 +\frac{3i}{2}B y^2)^T,\\
f_2(x) & = & (0, \frac{3}{4}B^2y)^T,\\
f_3(x) & = & (0, -\frac{i}{8}B^3)^T,
\end{eqnarray*}
\(f_{k}(x)\) equal to    the complex conjugate of \(f_{-k}(x)\) for \( k= -1, -2, -3 \),   and \( f_k(x)=0\)  for \( |k| > 3\).

\begin{figure}[t]
\begin{center}
\includegraphics[scale=.65]{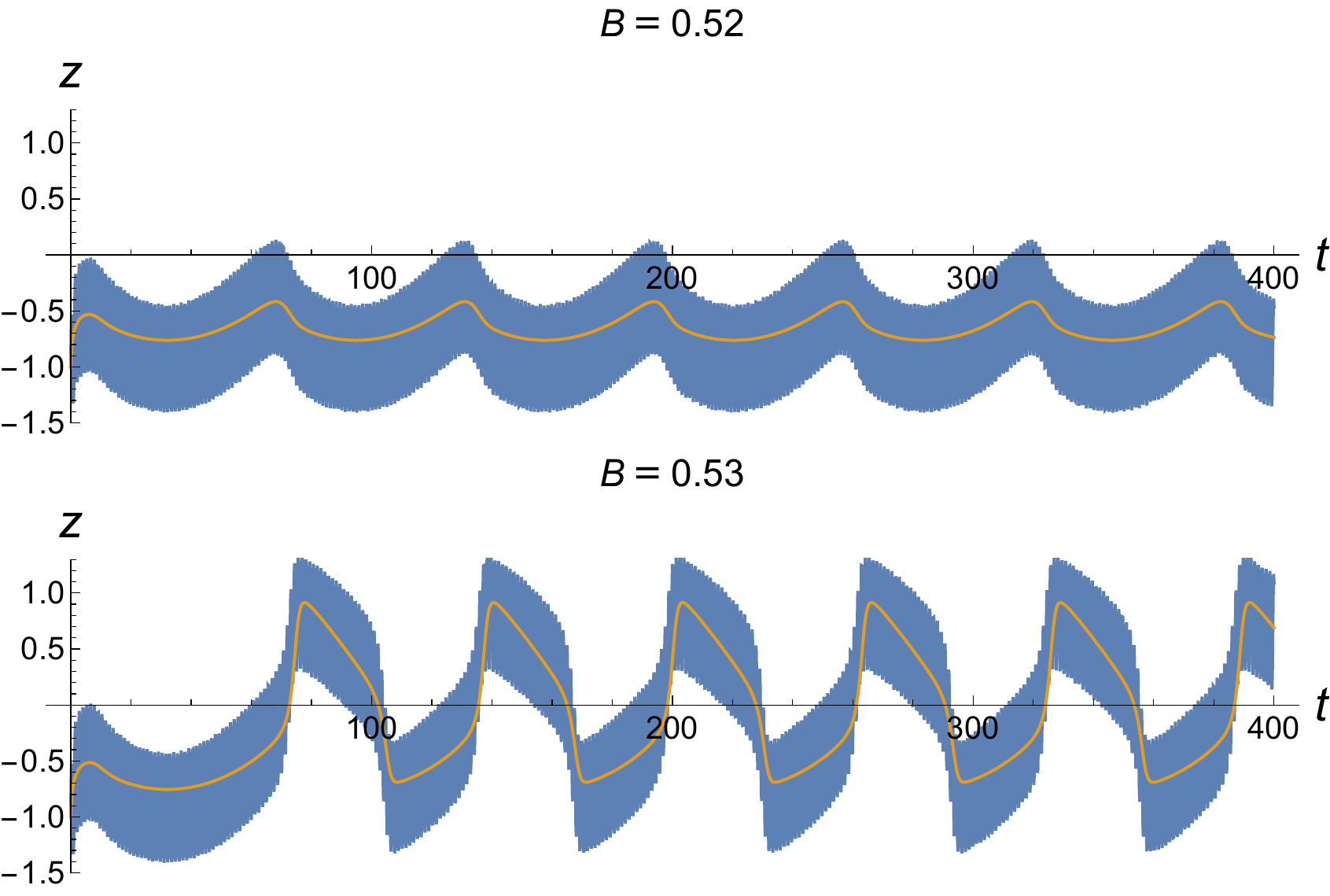}%
\end{center}
\caption{Vibrational resonance. Accurate numerical integrations of the original oscillatory differential equation  and of the averaged system based on words with \(\leq 2\) letters for two different values of the parameter \(B\) that measures the amplitude of the background vibration.  In each panel the oscillatory solution appears as a band due to its fast dynamics; the averaged solution (in the centre of the band) varies at a much lower rate. A small increase in \(B\)  from \(B=0.52\) (top panel) to \(B=0.53\) (bottom panel) lets the oscillator make substantially wider excursions without having to increment the amplitude \(A\) of the applied forcing.
}
\label{fig:fig}
\end{figure}

Truncating the averaged system \eqref{eq:edoaver} to keep only contributions corresponding to words of one letter results in
\[
\frac{d}{dt} \left[\begin{matrix}\Phi\\ Y \end{matrix}\right]=
\left[\begin{matrix}\nu\\ Y -\frac{3}{2} B^2 Y  -Y^3+A \cos \Phi \end{matrix}\right],
\]
i.e.\ the effect of averaging is only to remove the terms of \eqref{eq:vr2} that involve the background vibration. A simple computation of the required basis functions and coefficients that may be carried out by hand shows that truncating after two letter words yields \( (d/dt) \Phi=\nu\) (as expected) and
\begin{eqnarray}\nonumber
    \frac{d}{dt} Y &=& Y -\frac{3}{2} B^2 Y  -Y^3+A \cos \Phi \\
    \label{eq:twoletter}
     && {}+\frac{1}{\omega} \left(B-\frac{13}{6} B^3+B^5 -\frac{5}{2}B^3Y^2+3BY^4+6ABY \cos\Phi\right).
\end{eqnarray}
Thus the addition of two-letter words introduces corrections of size \(\mathcal{O}(1/\omega)\). The new {\em effective} potential
\[
-\frac{1}{2}Y^2 +\frac{3}{4} B^2 Y^2  +\frac{1}{4}Y^3+\frac{1}{\omega} \left(-BY+\frac{13}{6} B^3Y-B^5 Y +\frac{5}{6}B^3Y^3-\frac{3}{5}BY^5\right)
\]
has lost the \( Y\leftrightarrow -Y\) symmetry and the amplitude \( A\) of the slow driving forcing has been incremented by
an amount
\( 6ABY/\omega\).

Due to the absence of fast vibrations, \eqref{eq:twoletter} may be integrated numerically with high accuracy with a negligible computational effort. We have done so for
the values (taken from \cite{vr}) \(A = 0.2\), \(\nu = 0.1\), \(\omega = 5\). The initial condition is set to be
\(-1\), so that the oscillator starts at one of the minima of the potential well of \eqref{eq:vr1}. Figure~\ref{fig:fig}  shows accurate numerical solution of \eqref{eq:vr1} and \eqref{eq:twoletter} for two values of \( B\).   In each panel the oscillatory solution appears as a band due to its fast dynamics; the averaged solution varies at a much lower rate and is sufficient to describe the behaviour of the system.
When \(B=0.52\), the motion is essentially confined to the basin of the potential minimum at \(-1\); a small \(2\%\) increment in \(B\)  is sufficient to let the moving particle visit the basin of the minimum at \(1\) without having to increase the amplitude \( A\) of the driving forcing, thus demonstrating the existence of vibrational resonance. Note that, once the stationary regime is attained, the averaged solution is periodic with period \(2\pi/\nu = 20\pi\) and, as discussed above, does not have the \( Y\leftrightarrow -Y\) symmetry.

\begin{table}[t]
\begin{center}
\begin{tabular}{c|c|c|c}
\(n\) & \# \(n\)-letter words    & Error & Error \\
& with \(f_w\neq 0\)&\(B=0.52\)&\(B=0.53\)\\
\hline
1 & 7 & 0.241 & 0.431
\\
2 & 35 & 0.080 & 0.481
\\
3&217 & 0.026 & 0.251
\\
4& 1,407 & 0.018 & 0.036
\\
5 & 9,345 & 0.009 & 0.015
\\
6 & 62,951 & 0.004 & 0.008
\\
7&427,889 & 0.003 & 0.005
\end{tabular}
\end{center}
\caption{Higher-order averaging}
\label{table:table}
\end{table}

As pointed out before, the  solution \( Y\) of \eqref{eq:twoletter} approximates the true solution \(z\) at stroboscopic times. If approximations to \(z\) at non-stroboscopic times are also of interest, they may be easily obtained (without additional numerical integrations) by applying to the numerically computed \( Y\) the change of variables \eqref{eq:change} (see \eqref{eq:main}). If the change is truncated to exclude words with three or more letters, for the runs in Figure~\ref{fig:fig} the maximum on \(0\leq t\leq 400\) of the magnitude of the discrepancy between the true \(z\) and the approximation obtained in this way  is 0.080 for \(B= 0.52\) and \(0.481\) for \(B=0.53\).

Additional corrections of sizes \(\mathcal{O}(1/\omega^2)\), \(\mathcal{O}(1/\omega^3)\), \dots\ may be added to \eqref{eq:twoletter} by considering words with three, four, \dots\ letters. Explicit formulas for the
\(\mathcal{O}(1/\omega^2)\) terms (three-letter words) are given in \cite{kurusch}. The number of words to be considered increases exponentially with the length of the word. The second column of Table~\ref{table:table} gives, for each \(n=1, 2,\dots, 7\), the number of words  with \(n\) letters for which the associated basis function \(f_w\) is not identically zero. By using a computer algebra programme we have found the corrections of size \(\mathcal{O}(1/\omega^2)\), \dots, \(\mathcal{O}(1/\omega^6)\). The corresponding averaged systems are not reproduced here for obvious reasons.

For the parameter values used in Figure~\ref{fig:fig},  we have integrated numerically the averaged system based on words with \(\leq n\) letters \( n = 1, 2,\dots, 7\). Also listed in Table~\ref{table:table} is the maximum
error over \(0\leq t\leq 400\) when \( z \) is approximated by applying the change of variables including  words of 1, 2, \dots, \(n\) letters to the numerical solution of the averaged system employing the same words. It is apparent that the higher-order averaged systems, which, as emphasized before, may be easily integrated numerically, provide accurate approximations to the true oscillatory solution.

\bigskip
{\bf Acknowledgements.} We are thankful to M. A. F. Sanju\'an for bringing this problem to our attention and carefully reading the manuscript.
A. Murua and J.M.
Sanz-Serna have been supported by proj\-ects MTM2013-46553-C3-2-P and MTM2013-46553-C3-1-P from Ministerio de Eco\-nom\'{\i}a y Comercio, Spain. Additionally A. Murua has been partially supported by the Basque Government  (Consolidated Research Group IT649-13).

\end{document}